\documentclass[11pt]{article}

\usepackage{amsmath,amsthm,verbatim,amssymb,amsfonts,amscd, graphicx}
\usepackage{graphics}
\theoremstyle{plain}
\newtheorem{theorem}{Theorem}
\newtheorem{corollary}{Corollary}
\newtheorem{lemma}{Lemma}

\newtheorem*{surfacecor}{Corollary 1} 

\theoremstyle{definition}
\newtheorem{definition}{Definition}

\newcommand{\Z}{\mathbb{Z}}

\newcommand{\R}{\mathbb{R}}
\newcommand{\C}{\mathbb{C}}

\begin{document}

\title{The Thurston polytope for four-stranded pretzel links}
\author{Joan E. Licata}
\maketitle

\begin{abstract}
In this paper we use Heegaard Floer link homology to determine the dual Thurston polytope for pretzel links of the form $P(-2r_1-1, 2q_1, -2q_2, 2r_2+1), r_i, q_i \in {\Z}^+$.  We apply this result to determine the Thurston norms of spanning surfaces for the individual link components, and we explicitly construct norm-realizing examples of such surfaces.  
\end{abstract}

\section{Introduction}

The Thurston norm is a fundamental measure of the complexity of second homology classes in a three manifold.  When the manifold is a knot complement, $H_2(S^3, K)$ is spanned by any Seifert surface for $K$, and the $genus$ of $K$ is defined to be the least genus of any Seifert surface for $K$.  In the case of a link in the three-sphere, the analogous objects are surfaces whose boundaries lie on the components of the link.   In this paper we study the two-component pretzel links $P(-2r_1-1, 2q_1, -2q_2, 2r_2+1), p_i, q_i \in {\Z}^+$, which we denote by $P_{q_1,r_1,q_2,r_2}$, and we determine the dual Thurston polytopes of their complements.  


Let $B_T^*(P_{q_1,r_1,q_2,r_2})$ be the unit ball with respect to the dual Thurston norm on $H_1(S^3-P_{q_1,r_1,q_2,r_2} ;{\R})$.  We specify as a basis for $H_1(S^3-P_{q_1,r_1,q_2,r_2};{\R})$ the set of oriented meridians of the link components.  Identify the meridians with the coordinate vectors of ${\R}^2$ so that the meridian of the unknotted component is the vector $(1,0)$, and the meridian of the knotted component is the vector $(0,1)$.  Let $r_B$= max $(r_1, r_2)$, $r_S$= min $(r_1, r_2)$, $q_B$= max $(q_1, q_2)$, and $q_S$= min $(q_1, q_2)$.  In the special case $q_1=q_2=q$ and $r_1=r_2=r$, we denote $P_{q_1,r_1,q_2, r_2}$ by $P_{q,r}$.  

\begin{theorem}
$B_T^*(P_{q,r})$ is the convex hull of the following points: 
\[
\begin{array}{ll}
(-1, 4r-1) & (1, -4r+1) \\
(2q-3, 2q-3) & (-2q+3, -2q+3) \\
(2q-1, 2q-3) & (-2q+1, -2q+3) \\ 
(2q-1, -4r+2q-1) & (-2q+1, 4r-2q+1)
 \end{array} 
 \]
$B_T^*(P_{q_1, r_1, q_2, r_2})$ is the convex hull of the following points: 
\[
\begin{array}{ll}
(q_B-q_S-1, 2r_1+2r_2+q_B-q_S-1) & (-q_B+q_S+1, -2r_1-2r_2-q_B+q_S+1) \\
(q_1+q_2-3, 2r_B-2r_S+q_1+q_2-3) & (-q_1-q_2+3, -2r_B+2r_S-q_1-q_2+3) \\ 
(q_1+q_2-1, 2r_B-2r_S+q_1+q_2-3) & (-q_1-q_2+1, -2r_B+2r_S-q_1-q_2+3) \\
(q_1+q_2-1, -2r_1-2r_2+q_1+q_2-1) & (-q_1-q_2+1, 2r_1+2r_2-q_1-q_2+1)
 \end{array} 
 \]
 (If $2q$ (respectively, $q_1+q_2$)  $\le 2$, the points in the second line should be omitted.)
 \end{theorem}
 
\begin{figure}[h]
\begin{center}
\scalebox{.7}{\includegraphics{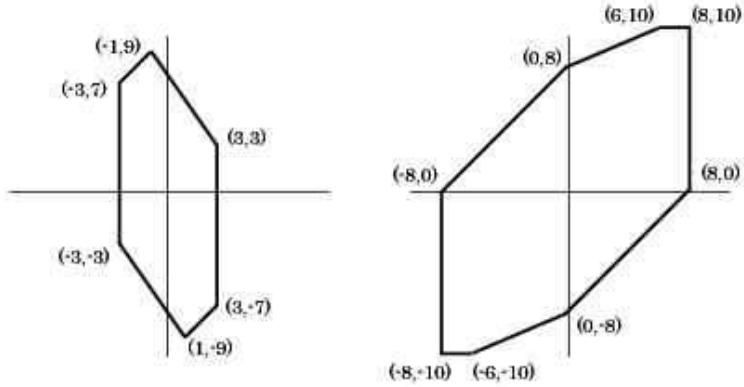}}  
 \end{center}
\caption{Left: $B_T^*(P_{2,2,2,3})$. Right: $B_T^*(P_{4,1,5,3})$.}\label {fig:examples}
\end{figure}

 Figure~\ref{fig:examples} shows two examples of dual Thurston polytopes.
 
McMullen (\cite{M}) showed that the dual Thurston polytope of the complement of a link $L$ contains the Newton polytope of the multivariable Alexander polynomial $\Delta(L)$.  The links $P_{q,r}$ have trivial Alexander polynomials (Corollary~\ref{cor:alexander}), so their Newton polytopes offer no information about $B_T^*(P_{q, r})$.  However, work of Ozsv\'{a}th, Szab\'{o}, and Ni has shown that link Floer homology detects the Thurston norm (\cite{OSz4}, \cite{N}), and we prove Theorem 1 by determining the filtration support for the Heegaard Floer link homology of $P_{q_1, r_1, q_2, r_2}$.

$B_T^*(P_{q_1, r_1, q_2, r_2})$ is the dual norm ball in $H_1(S^3-P_{q_1, r_1, q_2, r_2}; {\R})$ of the unit ball with respect to the Thurston norm in $H_2(S^3, P_{q_1, r_1, q_2, r_2}; {\R})$.  Thus, the Thurston norm of a homology class in $H_2(S^3, P_{q_1, r_1, q_2, r_2}; {\R})$ represented by a unit vector \textbf{u} is the half the length of the projection of $B_T^*(P_{q_1, r_1, q_2, r_2})$ onto a line parallel to \textbf{u}.  The theorem therefore implies the following:

\begin{corollary}\label{cor:surface}
The unknot component $U$ of $P_{q_1, r_1, q_2, r_2}$ bounds a surface $F_U$ with Euler characteristic $1-q_1-q_2$,  The knotted component $K$ of $P_{q_1, r_1, q_2, r_2}$ bounds a surface $F_K$ with Euler characteristic $- max \{ 2r_1+2r_2+q_B-q_S-1,  2r_B-2r_S+q_1+q_2-3\}.$  $F_U$ and $F_K$ have maximal Euler characteristic in their respective homology classes.
\end{corollary}

  The case $2r_1+2_2+q_B -q_S-1 > 2r_B-2r_S+q_1+q_2-3$ is particularly interesting, as $-\chi(F_K)$ then equals Thurston norm of the generator of $H_2(S^3, K)$.  This equality (proved in section 7) implies a minimal-complexity surface $(F_K, \partial F_K) \subset (S^3- P_{q,r}, \partial N(K))$ which is also a minimal-complexity Seifert surface for $K$ in $S^3$.  In particular, the links $P_{q, r}$ are split into two classes by the inequality $4r -1  \gtrless 2q-3$, and within each class the Thurston norm of $F_K$ is controlled exclusively by one of $q$ or $r$: provided that $2r+1$ is great enough, the Euler characteristic of $F_K$ does not depend on the linking with the unknotted component.  The opposite is true when $4r-1 < 2q-3$; then, the Euler characteristic of $F_K$  is a function solely of $q$.
  
  We begin with a discussion of the Thurston norm and pretzel links.  In section 4 we explicitly construct the norm-realizing surfaces $F_U$ and $F_K$.  Section 5 introduces the definitions and basic properties of Heegaard Floer link homology, and we prove Theorem 1 in the final section.

The author would like to thank Peter Ozsv\'{a}th for his encouragement and support.  In addition, the following individuals have been the source of many helpful conversations: John Baldwin, Elisenda Grigsby, Paul Melvin, Jiajun Wang, and Shaffiq Welji.

\vskip.5cm
\noindent

\section{The Thurston norm}

In \cite{T}, Thurston defined a semi-norm on the second homology of a three-manifold $M$ with boundary.  This norm measures the minimal complexity of a surface representing a fixed second homology class.  According to a result of Gabai (\cite{G}), it suffices to consider embedded surfaces.    In the case where $M$ is a link complement, we will take the "Thurston norm of $L$" to mean the Thurston norm of the manifold $S^3-L$.  

More precisely, let $S$ be an embedded (possibly disconnected) surface with components $s_i$, and define the $complexity$ of $S$ to be
\[\chi_-(S)=\sum_{ \{i \colon \chi(s_i) \le 0 \} } -\chi (s_i).\]

\begin{definition}
For $\sigma \in H_2(S^3, L)$, the Thurston norm of $\sigma$ is given by
\[||\sigma||_T=min\{ \chi_-(S) \colon [S]= \sigma \}.
\] 
\end{definition}
Thurston proved that the function $||\cdot||_T$ extends to a seminorm on $H_2(M,L; {\R})$ taking values in ${\R}$.  We specify a basis for $H_2(M,L)$ in which each element is a spanning surface for one component of the link and is disjoint from the other components.  If we identify each basis element with a unit vector on a coordinate axis of ${\R}^{|L|}$, the unit ball with respect to the Thurston norm is a polytope in ${\R}^{|L|}$.  

\begin{definition}
The dual Thurston norm $|| \tau ||_T^*$ for $\tau \in H^2(S^3, L; {\R}) \cong H_1(S^3-L; {\R})$ is given by
\[ || \tau ||_T^*= sup \{  |\tau (\sigma) | \}
\]
where the supremum is taken over $\sigma$ with $|| \sigma||_T =1$.
\end{definition}

The unit ball with respect to $|| \cdot ||_T^*$ is the dual polytope of the Thurston norm ball; we denote the dual Thurston norm ball by $B_T^*(L)$.

The Alexander polynomial $\Delta(L)$ of a link $L$  has one variable $x_i$ associated to the meridian $\mu_i$ of each component of $L$.  Sending $( \{ x_i \}, \cdot )$ to $( \{ \mu_i \}, + )$ identifies each monomial summand of $\Delta(L)$ with an element of $H_1(S^3-L)$.  The convex hull of these lattice points is known as the Newton polytope of $L$.  A result of McMullen [M] relates $B_T^*(L)$ to the Newton polytope, but in general $\Delta(L)$ does not determine $B_T^*(L)$.  Heegaard Floer link homology categorifies the multivariable Alexander polytonomial, and this stronger invariant determines the dual Thurston norm ball completely (\cite{OSz4}).

\section{Pretzel links}

Let $B_{a_i}$ be the two-strand braid $\sigma^{a_i}$, where $\sigma$ is the braid group generator and $a_i \in {\Z}$.  Adding bridges to cyclically connect the elements of a collection $\{ B_{a_i} \}^n_{i=1}$ gives the pretzel link $P(a_1, a_2...a_n)$. Figure 2 shows an example.

Ê
If the $a_i$ all have the same sign, then the link is alternating; in this case, the Alexander polynomial of the link determines its Thurston norm ([OSz4]).  For the remainder, we restrict our attention to $4$-tuples of twist coefficients whose signs alternate and parities do not. Although the four coefficients determine the link, there are equivalence classes of $4$-tuples that give isotopic links.  In particular, $P(2r_1+1, -2q_1, 2q_2, -2r_2-1)$ and $P(-2q_1, 2r_1+1, -2r_2-1, 2q_2)$ are isotopic, although $P(2r_1+1, -2q_1, 2q_2, -2r_2-1)$ and $P(-2q_1, 2r_1+1, 2q_2, -2r_2-1)$ are not.  

The Thurston norm of $L$ depends only on the homeomorphism type of the manifold and is therefore independent of the orientation of the link.   Since $B_T^*(L)$ is symmetric with respect to the origin, it is independent of the choice of orientation for the link components and their meridians.  Furthermore, the complements of a link and its mirror (the link $\bar{L}$ formed by changing every crossing in a planar projection of $L$) are homeomorphic.  Since the Thurston norm of a link and its mirror agree, it suffices to consider only $P(-2r_1-1, 2q_1, -2q_2, 2r_2+1)$ ($q_1, r_i \in {\Z}+$) which we denote by $P_{q_1,r_1, q_2, r_2}$.  This four-parameter family of links always has one knotted component, which we denote by $K$, and one unknotted component, which we denote by $U$.  For simplicity, we write $P_{q,r}$ for $P_{q,r,q,r}$.  

\section{Surfaces in the link complement}\label{sec:surf}   

 Surfaces representing basic classes in $H_2(S^3, P_{q_1,r_1, q_2, r_2})$ are Poincar\'{e} dual to the meridians of the various link components, and we refer to these as $spanning \  surfaces$.  Spanning surfaces are a generalization of Seifert surfaces, as the former may be punctured by other components of the link.  We present here a construction of the norm-minimizing spanning surfaces $F_U$ and $F_K$  in $S^3-P_{q_1,r_1, q_2, r_2}$.
 
Rephrasing Corollary~\ref{cor:surface}  in the language of Section 2, we have

\begin{surfacecor}  There exist surfaces $F_U$ and $F_K$ in $S^3-P_{q_1,r_1,q_2, r_2}$ satisfying the following:
\[-\chi(F_U)=||(1,0)||_T=q_1+q_2-1 \] \[-\chi(F_K)= ||(0,1)||_T=max \{ 2r_1+2r_2+q_B-q_S-1,  2r_B-2r_S+q_1+q_2-3\} .\]  
\end{surfacecor}

Any spanning surface for the unknotted component $U$ represents the homology class $(1,0)$, and   according to Corollary~\ref{cor:surface}, $||(1,0)||_T =q_1+q_2-1$.  The minimal-complexity surface $F_U$ is therefore a disc with $q_1+q_2$ punctures, and it is realized by the obvious Seifert surface for $U$ when $P_{q_1, r_1, q_2, r_2}$ is in the standard position.  See Figure 2 for an example.

\begin{figure}[!h]
\begin{center}
\scalebox{.8}{ \includegraphics{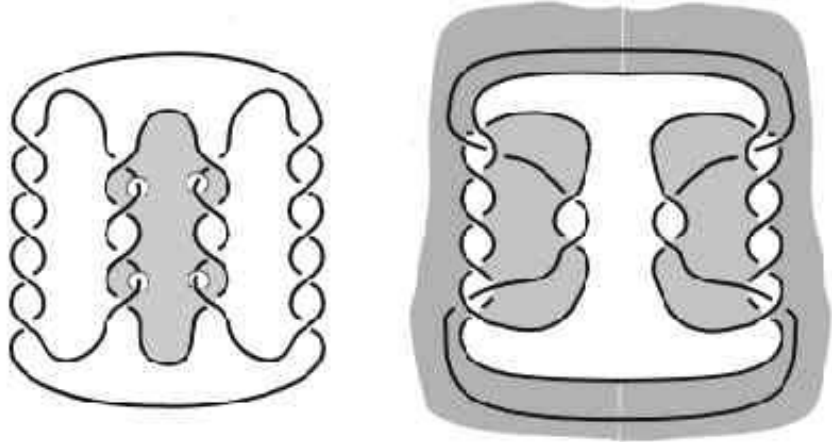}}
\end{center}
\label{fig:surface}
\caption{Left: $P_{2,2}$ is shown in standard position with a minimal-complexity spanning surface for the unknotted component.  The knotted component punctures the surface four times.  Right: The unknotted component is isotoped so as to be disjoint from a norm-realizing spanning surface for the knotted component.}
\end{figure}
Turning to $F_K$, we first consider minimal-complexity representatives of $(0,1) \in H_2(S^3, P_{q,r})$. As indicated in Corollary~\ref{cor:surface},  $||(0,1)||_T=max \{ 4r-1, 2q-3 \}$. The two cases correspond to the variations in the shape of $B_T^*(P_{q,r})$;  if $4r-1 < 2q-3$, the dual Thurston polytope is eight-sided, but if $4r-1  \ge 2q-3$, then $B_T^*(P_{q,r})$ is the convex hull of a proper subset of the points listed in Theorem 1.  Figure 2 illustrates a minimal spanning surface for the knotted component of $P_{2,2}$, but for larger values of $q$ and $r$, $F_K$ is best presented as a movie.   

If $P_{q,r}$ is presented in the standard projection, there is a natural Morse function $f: F_K \rightarrow [0,1]$ given by height on the page.  A Morse movie is a sequence of frames, where each frame is the preimage $f^{-1}(x)$ of some generic $x \in [0,1]$.  Any frame differs from the previous by isotopy or by a handle addition corresponding to a critical point of the Morse function, so the movie shows a descending sequence of horizontal slices through $F_K$ that captures the topology of the surface.   In our movie presentation of $F_K$, the initial frames will differ only by a isotopies dictated by the twisting of the strands.  Subsequently we will perform one-handle additions in order to arrive at a frame consisting of simple closed curves that may be capped off by two-handles.

The first frame, a slice taken immediately above the twisting, has three disjoint arcs parallel to the upper bridges of the link.  As the strands of the link twist, the leading edges spiral.  More precisely, consider the vertical strands of the link as a braid $B$ in ${\C} \times [0,1]$, and let $\mathbf{p} =B \cap ({\C} \times 0)$.  This braid defines a family of isotopies $B_s:( {\C}, \mathbf{p}) \rightarrow {\C} \times s, \ s \in [0,1]$ such that the image of $\mathbf{p}$ in ${\C} \times [0,1]$ is $B$. This family is well-defined up to isotopy fixing ${\C} \times \{0,1\}$.  The movie frame immediately below all the twisting is the image under $B_1$ of the "flat" arcs in the first frame.  This frame is clearly isotopic to the original one, but beginning in the next frame, we perform a sequence of saddle moves by adding one-handles.  

\begin{figure}
\begin{center}
\scalebox{.5}{ \includegraphics{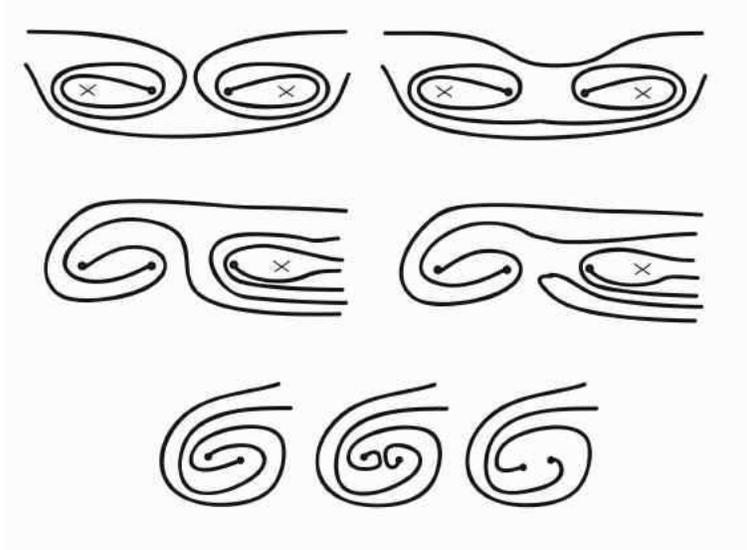}}
\end{center}
\caption{Top: $S_1$. Center: $S_2$ on the left side of the frame.  Bottom: $S_3$ on left side of the frame. Each strand of $U$ is represented by an x, and each strand of $K$ by a dot.}\label{fig:s1s2}
\end{figure}

\begin{figure}
\begin{center}
  \scalebox{.8}{\includegraphics{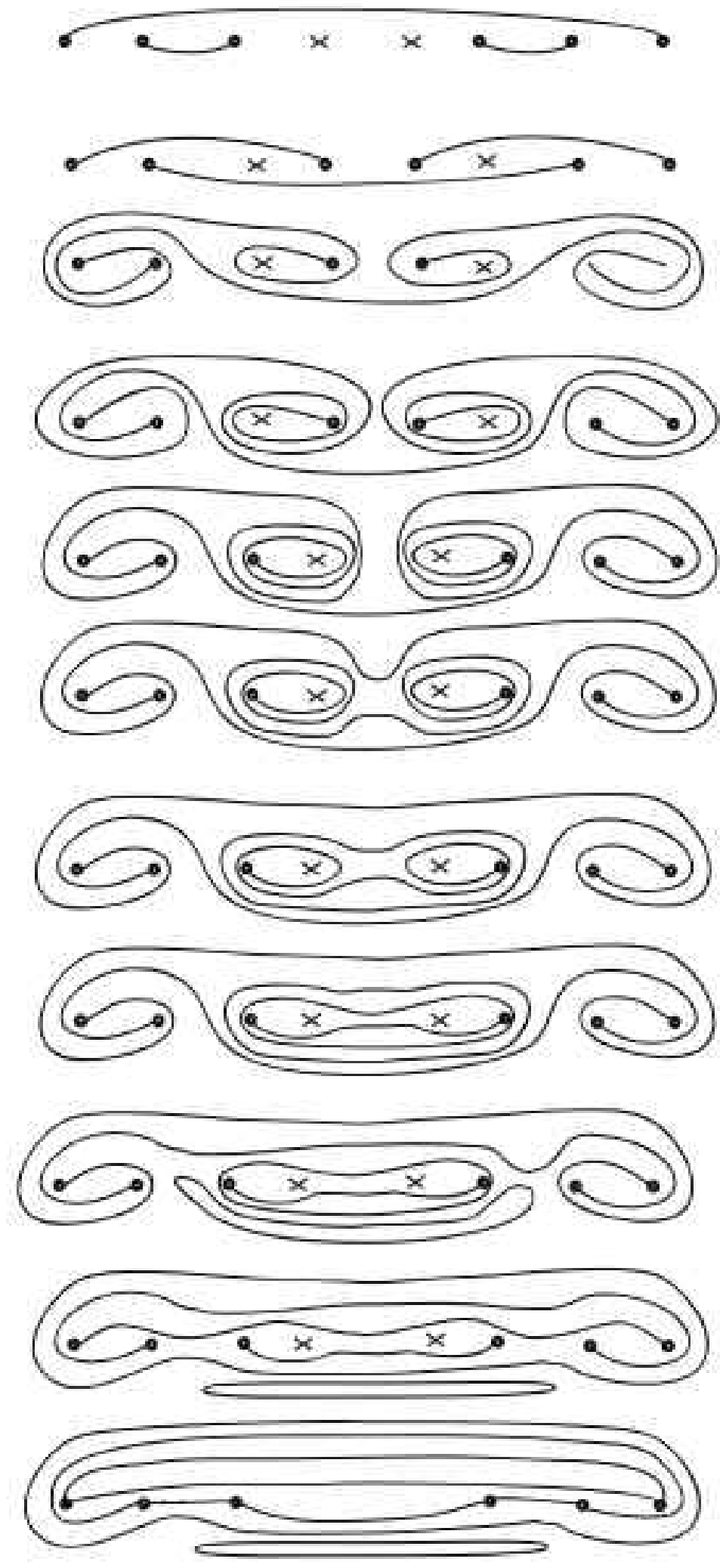}}  
  \end{center}
  \caption{The figure is a movie presentation of $F_K$ for $P_{3,1}$.  The x's are strands of $U$, and the dots are strands of $K$.  The surface has seven saddles and requires four death moves, so the Euler characteristic is -3.}\label{fig:movie2}
\end{figure}

Winding $K$ around $U$ (the braid action $\sigma_3^{2q}\sigma_5^{-2q}$) produced a pair of spirals in the center of the frame; define $S_1$ to be the saddle move applied at a pair of nearest points on these arcs.  Define $S_2$ as the saddle move that joins an arc from a lateral spiral to a nearest arc from the central spiral.  See Figure~\ref{fig:s1s2} for an illustration of these moves.   When $4r-1=2q-3$, performing $q \ S_1$ moves and $2(q-1) \ S_2$ moves yields a frame with $q-1$ simple closed curves and three arcs.  However, the arcs can be connected to form an additional two simple closed curves by extending the surface to the lower bridges of the link.  Capping off all the simple closed curves with two-handles gives a surface bounding $K$ with Euler characteristic $-2q+3$.  One can check that the surface constructed is orientable, and thus represents the desired homology class.  Figure~\ref{fig:movie2} shows an example of a movie presentation for $F_K$.  

In general, a third saddle move is necessary to fully resolve the diagram into simple closed curves.  $S_1$ produces a collection of concentric circles around an arc, and $S_2$ acts on the sides of these circles and arcs from the lateral spirals.  If $4r-1 > 2q-3$, however, the concentric circles will be exhausted before the lateral spirals, and we let $S_3$ be the move illustrated at the bottom of Figure~\ref{fig:s1s2} which reduces the lateral spiraling.  Applying $S_3$ to the remaining arcs of the lateral spirals yields a surface of Eular characterisitic $-2q+3 -2(2r-q+1)= -4r +1$, as desired.  

Finally, when $4r-1<2q-3$ we restrict the action of the braid group on the first frame to $\sigma_1^{-2r-1}\sigma_3^{4r+2}\sigma_5^{-4r-2}\sigma_7^{2r+1}$.  For each additional full twist of $U$ around $K$ ($\sigma_3^{2}$ or $\sigma_5^{-2}$), we introduce a puncture in the surface rather than a further isotopy of the leading edge.  Applying Morse moves as in the $4r-1=2q-3$ case then yields a Seifert surface for $K$ whose Euler characteristic is $1-4r-(2q-4r-2)=-2q+3$.

Although we described a construction for $P_{q,r}$, it extends to $P_{q_1, r_1, q_2, r_2}$ with the stipulation that equal powers of $\sigma_3$ and $\sigma_5^{-1}$ act on the first frame; any additional twists of $U$ around $K$ puncture the surface.  This power is $2q_S$ if $2r_S \ge q_S-1$, and the power is $4r_S+2$ if $2r_s < q_S-1$.  Maximal application of $S_1$ and $S_2$ (and possibly some $S_3$ moves) then yields an orientable surface with the desired complexity.  

\section{Heegaard Floer homology}\label{HF}
In order to prove Theorem 1, we will apply a result of Oszv\'{a}th and Szab\'{o} relating $B_T^*(L)$ to the Heegaard Floer homology of the link $L$.  We supply here a minimum of background on Heegaard Floer theory, and we refer the reader to \cite{OSz1}, \cite{OSz2},  \cite{OSz3}, and \cite{OSz4} for details.  


Let $\{ \alpha_i \}_{i=1}^{g+|L|-1}$ and $\{ \beta_i \}_{i=1}^{g+|L|-1}$ be collections of disjoint, transversely-intersecting,  simple closed curves on the genus-$g$ surface $\Sigma$.  Suppose that the $\{ \alpha_i \}$ and $\{ \beta_i \}$ each span $H_1(\Sigma)$, and that the Heegaard diagram induced by any pair of spanning subsets specifies $S^3$.  Let $\textbf{z}=\{ z_i \}$ and $\textbf{w}=\{ w_i \}$ be $|L|$-tuples of basepoints such that each component of $\Sigma - \{ \alpha_i \}$ (respectively, $\Sigma -\{ \beta_i \}$) contains a pair $(z_i, w_i)$.   Connect each ($z_i$, $w_i$) pair by two arcs, one in the complement of the $\{ \alpha_i \}$ and one in the complement of the $\{ \beta_i \}$.   Pushing these arcs into the corresponding handlebodies gives a simple closed curve $L_i$ in $S^3$.  If $L= \cup \{L_i\}$, we say that  $(\Sigma, \{ \alpha_i \}, \{ \beta_i \}, \textbf{z}, \textbf{w})$ is a $Heegaard\ diagram \ compatible \ with\ L$.

Given a space $X$, one constructs its $n^{th}$ symmetric product $Sym^n(X)$ from the product of $n$ copies of $X$ by modding out by the action of the symmetric group.  Define $\mathbb{T}_{\alpha}$ (respectively, $\mathbb{T}_{\beta}$) to be the image in $Sym^{g+|L|-1}(\Sigma)$ of the product $\alpha_1 \times \alpha_2 .... \times \alpha_{g+|L|-1}$ ($\beta_1 \times \beta_2 .... \times \beta_{g+|L|-1}$).  Since the $\alpha$ and $\beta$ curves intersect transversely in $\Sigma$, the associated tori intersect transversely in $Sym^{g+|L|-1}(\Sigma)$.  Let $\widehat{CFL}(L)$ be the free abelian group generated by the intersection points $\mathbb{T}_{\alpha}\cap\mathbb{T}_{\beta}$.  We will equip $\widehat{CFL}(L)$ with a boundary map to give it the structure of a chain complex, but we defer this construction until after a discussion of the multi-grading.

Supose that \textbf{x} and \textbf{y} are intersection points of $\mathbb{T}_{\alpha}\cap\mathbb{T}_{\beta}$, and let $\gamma_+$ (respectively, $\gamma_-$) be the arc of the unit circle with positive (negative) real component.  We denote by $\pi_2(\textbf{x}, \textbf{y})$ the set of all homotopy classes of maps from the unit disc to $Sym^{g+|L|-1}(\Sigma)$ that satisfy the following: 
\[
\left\{
   \begin{array}{@{}c@{}}
      \phi(i)= \textbf{x}\\
      \phi(-i)=\textbf{y}\\
      \phi(\gamma_+) \subset \mathbb{T}_{\beta}\\
      \phi(\gamma_-) \subset \mathbb{T}_{\alpha}
  \end{array}
 \right\}  
 \]

If $\phi$ is such a map,  we denote its image in the symmetric product by $\phi$ as well.
For a point $p$ in $\Sigma-\{\alpha_i\}-\{\beta_i\}$, the intersection number of $\phi$ with $p$ is given by
\[
n_p(\phi) = \#( \phi \cap (\{p\} \times Sym^{{g+|L|-2}}(\Sigma))),
\] 
where the intersection takes place in $Sym^{g+|L|-1}(\Sigma)$.

For ease of use, $\Sigma$ is preferable to $Sym^{{g+|L|-1}}(\Sigma)$, and we interpret the preceding definitions accordingly.  Thus, we note the existence of a bijection between generators of $\widehat{CFL}(L)$ and unordered $(g+|L|-1)$-tuple of intersection points on $\Sigma$ with the property that each of the $\alpha_i$ and $\beta_j$ circles contributes to exactly one intersection point.  The boundary of $\phi$ corresponds to a collection of arcs in the $\alpha$ and $\beta$ curves whose endpoints are the intersection points.  

To a disc $\phi$ in  $Sym^{g+|L|-1}(\Sigma)$ we may associate the immersed surface $D(\phi)$ in $\Sigma$.
\begin{definition}
Let $\Omega_i$ be the components of $\Sigma-\{\alpha_i\}-\{\beta_i\}$, and let $p_i \in \Omega_i$ be a point in each such component.  For a disc $\phi$, the domain of $\phi$ is the formal object
\[D(\phi)=\sum_{i=1}^N n_{p_i} \Omega_i
\]
\end{definition}

Given two $(g+|L|-1)$-tuples of intersection points in $\Sigma$ which represent generators \textbf{x} and \textbf{y}, one builds $D(\phi)$ for $\phi \in \pi_2( \textbf{x}, \textbf{y})$ by first tracing out alternating $\alpha$ and $\beta$ arcs connecting the intersection points.  Intersections corresponding to one of the generators are the only corner points allowed, so any other boundary components must be complete $\alpha$ or $\beta$ circles.  This boundary determines the multiplicities of the $\Omega_i$, and this method builds the domain for an arbitrary $\phi$ connecting \textbf{x} and \textbf{y}.  

Two domains connecting the same generators may differ by a copy of $\Sigma$ or by a connected component of $\Sigma-\{ \alpha_i \}$ or $ \Sigma - \{ \beta_i \}$.  A periodic domain is a homological relation among the $\alpha$ and $\beta$ circles; concretely, it is a collection of $\Omega_i$ whose boundary is a subset of $\{ \alpha_i \cup \beta_i \}$ such that no $\Omega_i$ containing a $z_i$ or $w_i$ is included.  If $\phi, \phi' \in \pi_2(\textbf{x}, \textbf{y})$ have the same multiplicities at all the basepoints, then $D(\phi)$ and $D(\phi')$ differ by some periodic domain.  

We are now in a position to define the multi-grading on $\widehat{CFL}(L)$.
\begin{definition}
Let $\phi$ be a disc connecting \textbf{x} to \textbf{y}, and let $F(\textbf{x}, \textbf{y})$ be the vector given by
\begin{equation} F(\textbf{x}, \textbf{y}) = (n_{z_1}(\phi)-n_{w_1}(\phi), n_{z_2}(\phi)-n_{w_2}(\phi)...n_{z_{|L|}}(\phi)-n_{w_{|L|}}(\phi)).  
\end{equation}
We say that $F(\textbf{x}, \textbf{y})$ is the filtration level of \textbf{x} relative to \textbf{y}.   This vector is independent of the choice of $\phi \in \pi_2(\textbf{x}, \textbf{y})$.
\end{definition}
$F$ partitions the generators of $\widehat{CFL}(L)$ into relative filtration levels: $\textbf{x} \backsim \textbf{y}$ if $F(\textbf{x}, \textbf{y})=0$.  There is a canonical identification of the filtration partition classes in $\widehat{CFL}(L)$ with ${\Z}^{|L|}$.  The filtration support of link homology is naturally symmetric around a center point, and we identify this with the origin.  With these coordinates, we have the symmetry relation
\begin{equation}\label{eq:sym}
 \widehat{HFL}(L, \textbf{h}) \cong  \widehat{HFL}(L, \textbf{-h}).
\end{equation}

The meridians of the link components generate $H_1(S^3-L)$, and in fact the quantity $n_{z_i}(\phi)-n_{w_i}(\phi)$ may be interpreted as the linking number of the $i^{th}$ meridian with a loop in $S^3-L$ associated to the generators \textbf{x} and \textbf{y}.  This in turn allows us to identify a filtration level \textbf{h} with an element of $H_1(S^3-L)$.   For \textbf{x} $\in \widehat{CFL}(L, \textbf{h})$ we refer to the $H_1(S^3-L)$ grading \textbf{h} as the filtration level of \textbf{x} in order to distinguish this from the Maslov, or homological, grading.  Although defined analytically, the relative Maslov grading of a pair of generators can be calculated via a combinatorial formulation presented first in \cite{L}.

Let $\textbf{x}=(x_1, x_2...x_n)$ and $\textbf{y} =(y_1, y_2...y_n)$ be generators of $\widehat{CFL}(L)$, and let $\phi$ be any element of  $\pi_2 (\textbf{x}, \textbf{y})$.  Consider $D(\phi)$ as a surface in $\Sigma$ with boundary in the four-valent graph formed by the $\alpha$ and $\beta$ curves.  Let $k$ be the number of corner points where $D(\phi)$ fills a single quadrant, and let $l$ be the number of corners where it fills three quadrants.  Furthermore, let $\bar{n}_{x_i}$ be the local multiplicity of the intersection point $x_i$; for example, this is $\frac{1}{4}$ if $x_i$ is a $k$-type corner point, and it is $n_{x_i}$ if $x_i$ is in the interior of $\Omega_i$.  The following formula computes the Maslov index $\mu(\phi)$ of $\phi$, where $\chi$ refers to the Euler characteristic of the surface:
\begin{equation}\label{eq:maslov}
	\mu(\phi)=\chi(D(\phi))-\frac{k}{4}+\frac{l}{4}+\sum_i \bar{n}_{x_i}+\sum_i \bar{n}_{y_i}-2\sum_i n_{w_i}(D(\phi))).
\end{equation}	
\begin{definition}
The (homological) grading of \textbf{x} relative to \textbf{y}  is the Maslov index of any disc $\phi \in \pi_2(\textbf{x}, \textbf{y})$.
\end{definition}

A complex structure on $\Sigma$ induces a complex structure on $Sym^{g+|L|-1}(\Sigma)$, and the Maslov index gives the expected dimension of the moduli space $\mathcal{M}(\phi)$ of holomorphic representatives of the homotopy class $\phi$.  Identifying the unit disc conformally with a vertical strip in the complex plane allows one to mod out by the ${\R}$ action that shifts the strip vertically.  We denote the resulting moduli space by $\frac{\mathcal{M}(\phi)}{{\R}}$. If $\mu(\phi)=1$, $\frac{\mathcal{M}(\phi)}{{\R}}$ is zero-dimensional.  (For a more detailed discussion of the moduli spaces, we refer the reader to \cite{OSz1}.)  

We are finally in a position to define the differential $\widehat{\partial} \colon \widehat{CFL}(L) 
\rightarrow \widehat{CFL}(L) $.  
\begin{definition} The boundary map $\hat{\partial}$ on  $\widehat{CFL}(L)$ is given by
\[
\hat{\partial} (\textbf{x}) =  \sum_{\textbf{y}} \sum_{\phi \in \pi_2(\textbf{x},\textbf{y}) \colon \mu(\phi)=1\  n_{\textbf{w}}(\phi)=n_{\textbf{z}}(\phi)=0} \#(\frac{\mathcal{M}(\phi)}{{\R}}) \textbf{y}
\]
\end{definition}

This map satisfies $\hat{\partial}^2=0$, so $(\widehat{CFL}, \hat{\partial} )$ is a chain complex whose homology we denote by $\widehat{HFL}(L)$.  Note that the boundary preserves the filtration level and lowers the grading by one, so in fact 
\begin{equation}
\widehat{HFL}(L) = \bigoplus_{\textbf{h} \in H_1(S^3-L, {\Z})}\widehat{HFL}(L, \textbf{h}).
\end{equation}

Recall from section 2 that the dual Thurston norm is defined on $H_2(S^3, L)^* \cong H_1(S^3-L)$.  This allows us to compare the dual norm ball to the set $\{ \textbf{h} \colon \widehat{HFL}(L, \textbf{h}) \ne 0 \}$).
  \begin{theorem}[OSz4]\label {thm:polytope}
The setwise sum in ${\Z}^{|L|}$ of the dual Thurston polytope and the cube of edge-length two is twice the convex hull of the filtration support of $\widehat{HFL}(L)$.  
\end{theorem}

\section{Proof of Theorem 1}

We begin by computing the dual Thurston polytope for the knotted component $K$ of $P_{q_1, r_1, q_2, r_2}$.  In addition to serving as a toy calculation illustrating the techniques, the result will prove useful in the sequent.  

\begin{figure}[h]
\begin{center}Ê 
\scalebox{.5}{\includegraphics{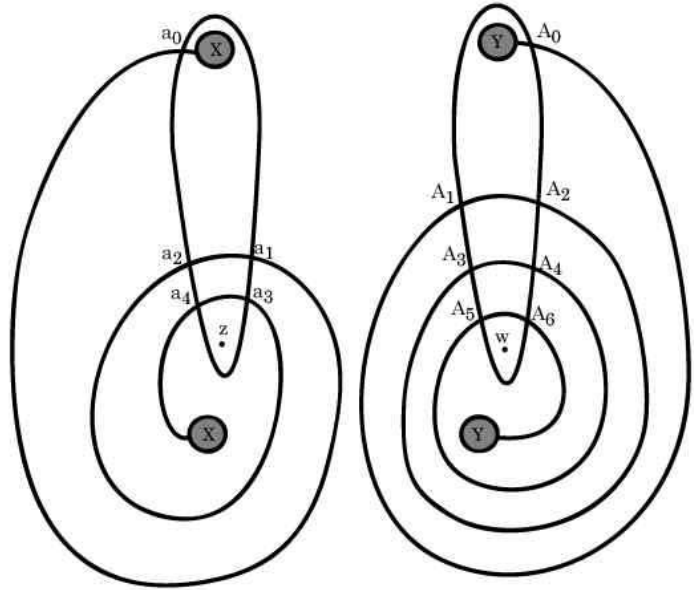}}
\end{center}
\caption{A genus-two Heegaard diagram compatible with $T(-2, 5) \# T(2,7)$.  The vertical pairs of shaded circles should be seen as the gluing discs for one-handles along which the $\beta$ curves run.}\label {fig:KHD}
\end{figure}

\begin{figure}
\begin{center}
\scalebox{.8}{\includegraphics{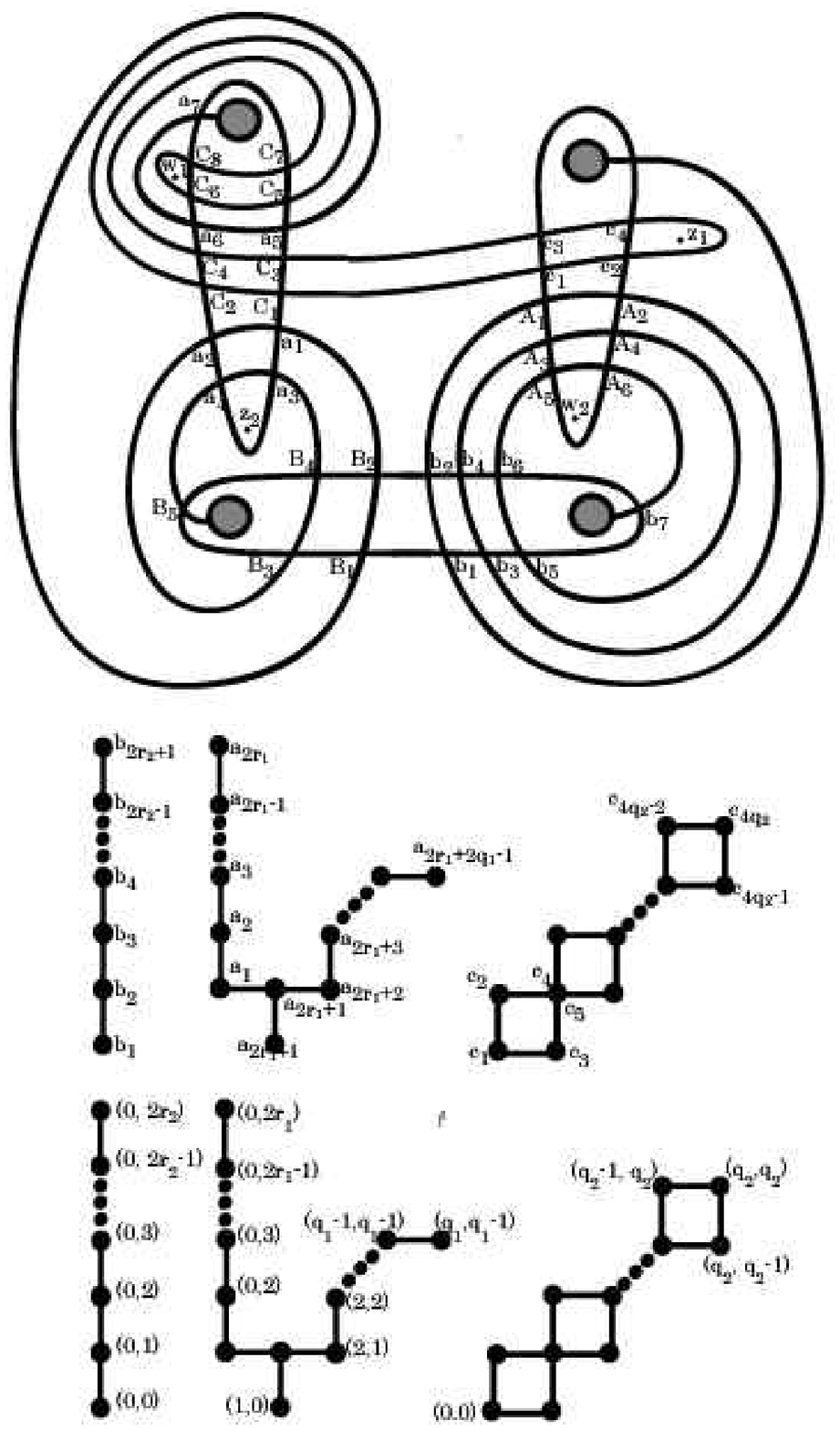}}
\end{center}
\caption{Top: Heegaard diagram for $P_{2,2,1,3}$. Bottom: Relative filtration data and temporary coordinates for lowercase intersection points.  The "filtration level" of a generator is the vector in ${\Z}^2$ given by summing the coordinates of its three constituent intersection points. Switching the subscripts 1 and 2 gives the filtration data for the uppercase generators.}\label{fig:P2213}
\end{figure}

$K$ is the connected sum of the torus knots $T(-2,2r_1+1)$ and $T(2, 2r_2+1)$.  Figure~\ref{fig:KHD} shows a Heegaard diagram compatible with $K$, and we see that each generator of $\widehat{CFL}(K)$ corresponds to a pair of intersection points $a_i \in \alpha_1 \cap \beta_1$ and $A_j \in \alpha_2 \cap \beta_2$.  Holding the $\alpha_2 \cap \beta_2$ intersection point constant, one determines the relative filtration levels of the $a_i$ by finding a domain connecting any pair.  Relative filtration is additive under concatenation of discs, so it suffices to show that successively numbered intersection points have relative filtration one.  The same holds for the points of $\alpha_2 \cap \beta_2$, so the filtration level of  the generator $(a_{2r_1}, A_{2r_2})$ is $2r_1+2r_2$ greater than that of $(a_{0}, A_{0})$.  These two generators are unique in their respective filtration levels, so each represents a non-trivial homology class.  The Heegaard Floer polytope, then, is the interval $[ -r_1-r_2, r_1+r_2]$.  (Recall the symmetry relation of Equation~\ref{eq:sym}.)  To translate this into information about the Thurston norm of $S^3-K$, we apply the result of Ozsv\~ath and Sz\~abo stated in Theorem~\ref{thm:polytope}; doubling the interval and subtracting 1, we find that the minimal complexity of any Seifert surface for $K$ is $2r_1+2r_2-1$.  If $K$ is drawn as the pretzel knot $P(2r_1+1, 0, -2r_2-1)$, Seifert's algorithm gives a surface realizing the Thurston norm.
Turning now to the proof of the main theorem, we use a Heegaard diagram of the type indicated in Figure~\ref{fig:P2213} to bound the Heegaard Floer polytope of $P(-2r_1-1, 2q_1, -2q_2, 2r_2+1)$ .  When $r_1 \ne r_2$, these bounds determine $B_T^*(P_{q_1, r_1, q_2, r_2})$, but we will ultimately need to employ a second Heegaard diagram to complete the calculation for $P_{q,r}$.

Generators of $\widehat{CFL}(S^3, P_{q_1, r_1, q_2, r_2})$ correspond to triples of intersection points on the diagram, and these split naturally into two types.  Denote triples of the form $(\alpha_1 \cap \beta_1, \alpha_2 \cap \beta_3, \alpha_3 \cap \beta_2)$ by lowercase letters  $(a_*, b_*, c_*)$, and triples of the form $(\alpha_1 \cap \beta_2, \alpha_2 \cap \beta_1, \alpha_3 \cap \beta_3)$ by uppercase letters  $(A_*, B_*, C_*)$. 
\begin{lemma}\label{filtrpr} 
The filtration levels of $(A_i, B_j, C_k)$ and $(a_i, b_j, c_k)$ are equal.
\end{lemma}

\begin{proof}[Proof of Lemma~\ref{filtrpr}]
Equation~\ref{eq:sym} determines the absolue filtration on $\widehat{HFL}(L)$, but until the homology is known, we have information only about the relative filtration.  For convenience, therefore, we adopt pro tem the coordinates indicated in Figure~\ref{fig:P2213}.  Specifically, we situate the graphs in ${\Z}^2$ by placing $b_1$ and $c_1$ at the origin, and $a_1$ at $(0,1)$.   The filtration level of a generator is given by the sum of the filtration levels corresponding to the three constituent intersection points, as  relative filtration is additive under concatenation of homotopy classes of discs. We can therefore say that the generator $(a_1, b_1, c_1)$ has filtration level $(0,1)$, while $(a_2, b_1, c_1)$ has filtration level $(0,2)$; as we see by comparing their $a_i$ intersections, this is consistent with the two generators' relative filtration being $(0,1)$.  

  There is a hexagon disjoint from all the basepoints which connects $(A_1, B_1, C_1)$ to $(a_1, b_1, c_1)$, so these generators are in the same filtration level.  Now suppose that $\phi$ is a domain connecting $(a_1, b_1, c_1)$ to another lowercase generator $(a_i, b_j, c_k)$ with filtration $F(\phi) = (n_{z_1}(\phi)-n_{w-1}(\phi), n_{z_2}(\phi)-n_{w-2}(\phi)...n_{z_n}(\phi)-n_{w-n}(\phi)$.  The mirror image of $\phi$ is a domain $\phi'$connecting $(A_i, B_j, C_k)$ to $(A_1, B_1, C_1)$ and satisfying $F(\phi')=-F(\phi)$.  Consequently each uppercase generator has the same filtration level as its lowercase counterpart.  Note that if $q_1 \ne q_2$ or $ r_1\ne r_2$, some generators will not be part of such pairs.
\end{proof}

\begin{corollary}\label{cor:alexander}
The multivariable Alexander polynomial vanishes for the links $P_{q,r}$.
\end{corollary}

\begin{proof}[Proof of Corollary~\ref{cor:alexander}] 
The graded Euler characteristic of Heegaard Floer link homology is the multivariable Alexander polynomial (\cite{OSz3}).  According to Lemma~\ref{filtrpr}, the generators of $\widehat{CFK}(S^3, P_{q,r})$ appear in pairs with the same filtration levels.  We claim that the gradings of such a pair of elements are always of opposite parity.  In this case, an application of the Euler-Poincar\'{e} principle shows that the graded Euler characteristic of the homology vanishes. 

To prove the claim, we first note that it holds for the pair $(a_1, b_1, c_1)$ and $(A_1, B_1, C_1)$, since the hexagonal domain cited above has Maslov index one.  Furthermore, if  $\phi$ and $\phi'$ are the mirror-image domains from the proof of Lemma~\ref{filtrpr}, then their Maslov indices are equal mod two; the only difference comes from the even-integer-valued term $2n_{w_i}$.  
\end{proof}

\begin{lemma}\label{lem:nonvan} For $(x,y) \in \{(0, 2r_1+2r_2+1), (q_B-1, 2r_1+2r_2+q_B), and (q_B, 2r_1+2r_2+q_B)\}$, the Heegaard Floer groups $\widehat{HFL}(S^3, P_{q_1,r_1, q_2, r_2}, (x,y)) $are nontrivial.
\end{lemma}

\begin{figure}[h]
\begin{center}
\scalebox{.8}{\includegraphics{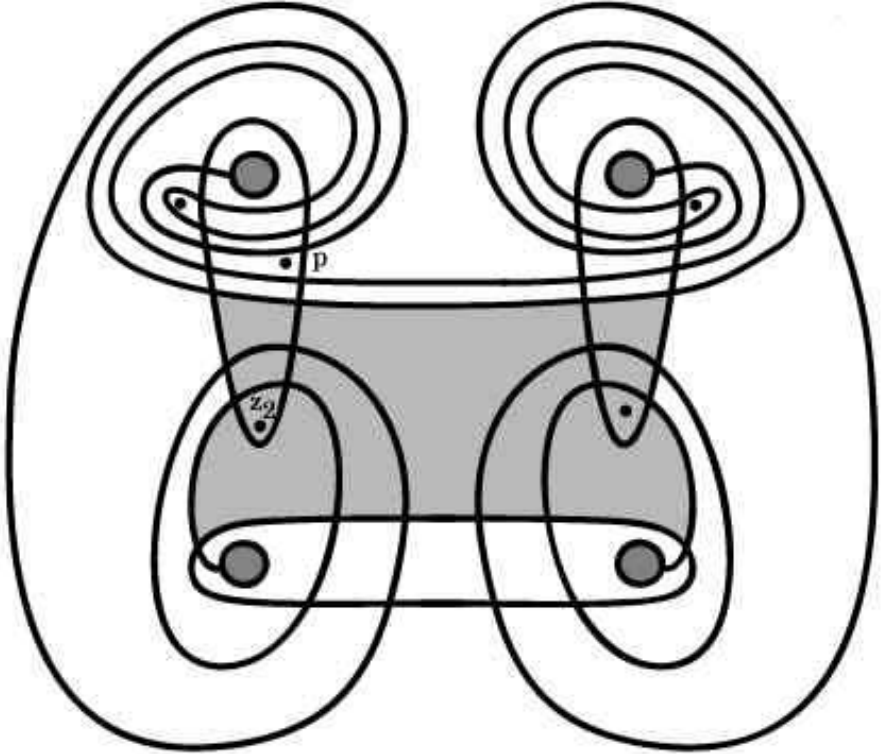}}
\end{center}
\label{HDEx}
\caption{The figure shows $D(\phi)$ for a disc connecting the generator $(A_{2r_2}, B_{2r_1+1}, C_2)$ to $(a_{2r_1}, b_{2r_2+1}, c_2)$.  Since $n_{z_2}(\phi) > n_p(\phi)$, $\widehat{HFL}(L, (0, 2r_1+2r_2+1) \ne 0$.}\label {fig:HDEx}
\end{figure}

\begin{proof}[Proof of Lemma~\ref{lem:nonvan}]
If $q_1  \ne q_2$, the filtration levels $(q_B-1, 2r_1+2r_2+q_B)$ and $(q_B, 2r_1+2r_2+q_B)$ each have a unique generator, so the homology is also one-dimensional.  If $q_1 = q_2$, there is a single pair of generators at each of these filtration levels and also at $(0, 2r_1+2r_2+1)$.  
We show $\widehat{HFL}(S^3, P_{q_1,r_1, q_2, r_2}, (x,y)) \ne 0$ by showing that for each of these pairs, no disc $\phi$ connecting them can have a holomorphic representative.  By definition, a boundary disc satisfies $n_{z_i}(\phi)=n_{w_i}(\phi)=0$, and a disc with a holomorphic representative must also satisfy $n_{p_i}(\phi) \ge 0$ for all $p_i$ in the complement of the $\alpha$ and $\beta$ circles.  ([OSZ1], Lemma 3.2)) For each of the generator pairs in question, one may show that there is no domain satisfying both these conditions.  

We demonstrate this calculation explicitly for $(x,y)=(0, 2r_1+2r_2+1)$; the other two proceed similarly.  We first consider a domain for an arbitrary $\phi_1 \in \pi_2((A_{2r_2}, B_{2r_1+1}, C_2), (a_{2r_1}, b_{2r_2+1}, c_2))$.  Figure~\ref{fig:HDEx} illustrates the domain on the diagram for $P_{2,2}$.  Note that $\phi_1$ has Maslov index one, so it potentially represents a map from the higher-graded capital generator to the lower-graded lowercase one. The $\alpha$ curves separate the surface into two connected components, as do the $\beta$ curves,  and these subsurfaces may be added or subtracted from $\phi_1$ to change the basepoint multiplicities.  Subtracting a copy of the component of $\Sigma- \{ \alpha_i \}$ containing the handles gives a new domain $\phi_2$ satisfying $n_{w_1}(\phi_2)=n_{z_1}(\phi_2)=n_{w_2}(\phi_2)=n_{z_2}(\phi_2)=0$.  As noted in Section~\ref{HF}, another domain with the same intersection numbers at the basepoints can differ from $\phi_2$ only by periodic domains.  On this diagram, the space of periodic domains is one-dimensional over ${\Z}$.  We note that the generating periodic domain has multiplicity zero at the extra marked point $p$, so adding a periodic domain cannot change $n_p (\phi_2)$.   Since $n_p (\phi_2) < 0$, the positivity  and null-intersection conditions cannot be simultaneously realized.  The moduli space $\frac{\mathcal{M}(\phi)}{{\R}}$ is therefore empty, and $\widehat{HFL}(S^3, P_{q_1,r_1, q_2, r_2}, (0, 2r_1+2r_2+1))\ne 0 $.  

\end{proof}         

\begin{proof}[Proof of Theorem 1]
     According to Theorem~\ref{thm:polytope} (\cite{OSz4}), $B_T^* (P_{q_1,r_1, q_2, r_2})$ is determined by the convex hull $\mathfrak{H}(\widehat{HFL})$ of the lattice points $(x,y)$ where $\widehat{HFL}(S^3, P_{q_1,r_1, q_2, r_2}, (x,y)) \ne 0$.  (Note that $(x,y)$ refers to a relative filtration level, as in the proof of Lemma~\ref{filtrpr}.)  We first consider $S \colon = \{(x,y) \colon \widehat{CFL}(S^3, P_{q_1,r_1, q_2, r_2}, (x,y)) \ne \emptyset) \}$ and its convex hull $\mathfrak{H}(S)$.   

     In an ideal world, the filtration support of the nontrivial homology might be equal to the filtration support of $\widehat{CFK}(P_{q_1,r_1, q_2, r_2})$ and the polytopes $\mathfrak{H}(S)$ and $\mathfrak{H}(\widehat{HFL})$ would coincide .  This is obviously false, however, as the convex hull of the filtration support is not symmetric.  Lowering our expectations a bit, we can ask if the two polytopes share enough of the same vertices to ensure that the symmetry relation determines the rest of of $\mathfrak{H}(\widehat{HFL})$, as indicated in Figure~\ref{fig:parttopes}.  In fact, when $r_1 \ne r_2$, this is almost the case.
     
  \begin{figure}
Ê \includegraphics[width=\textwidth]{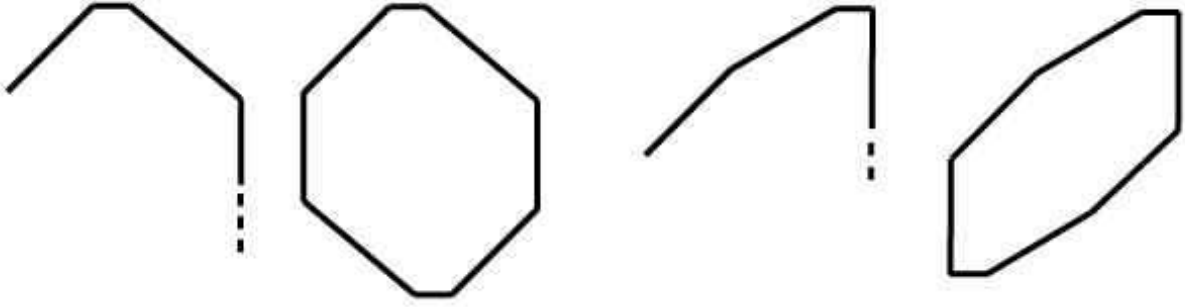}
\caption{Since the Heegaard Floer polytope is symmetric, half the sides determine the entire object.  The sides with rational slope agree with those of $\mathfrak{H}(S)$, although the vertical edge does not.  The pair of figures on the left correspond to the case $2r_S > q_S-1$, and the pair of figures on the right to the case $2r_S<q_S-1$. }\label{fig:parttopes}
\end{figure}
     
     The filtration data determine the extreme points of  $\mathfrak{H}(S)$.  For the minimal $x$ value, the point $(0, 2r_1+2r_2+1)$ maximizes $y$.  Varying the $c_i$ intersection shows that the points $(k, 2r_1+2r_2+1+k)$ are in $S$ for $0 \le k\le q_B-1$, and for a fixed $x$, these are the maximal such $y$ values.  Note also that $(q_B, 2r_1+2r_2+q_B)$ is realized by at least one generator, or by two if $q_1=q_2$. According to Lemma~\ref{lem:nonvan}, each of these points is also in $\mathfrak{H}(\widehat{HFL})$.  
     Moving right, 
\begin{equation}\label{eq:S}
max \{ y \colon (x,y) \in S, q_B+1\le x\le q_1+q_2 \}= q_1+q_2+2r_B-1.
\end{equation}
When $r_1 \ne r_2$, there is a unique generator in each of the filtration levels $(q_1+q_2-2, q_1+q_2+2r_B-1)$ and $(q_1+q_2, q_1+q_2+2r_B-1)$.  Since these filtration levels consequently have nontrivial homology, the top faces of $\mathfrak{H}(S)$ and  $\mathfrak{H}(\widehat{HFL})$ agree. 

To completely determine $\mathfrak{H}(\widehat{HFL})$ we need also the length of one of the vertical sides, but here the two polytopes differ: the two generators in filtration level $(q_1+q_2, q_1+q_2-2)$ are connected by a topological disc domain that has a holomorphic representative, so the homology is trivial.  However, we will show $\widehat{HFL}(P_{q_1,r_1,q_2, r_2}, (q_1+q_2, 2q_1+q_2-1)$ is nontrivial, proving that the vertical edges of $\mathfrak{H}(\widehat{HFL})$ have length $2r_B$.
     
The complexity of any spanning surface for $K$ in $H_2(S^3, P_{q_1,r_1, q_2, q_2})$ is bounded from below by the Thurston norm of a Seifert surface for $K$ in $S^3$.  Using the calculation from the beginning of the section, this means $||(0,1)||_T \ge 2r_1+2r_2-1$.  The total height of the Heegaard Floer polytope must therefore be at least $(2r_1+2r_2-1)+1=2r_1+2r_2$.  For the moment setting $q_1=q_2 =q\le 2r_S+1$ and computing the vertical component of the sloped edges, we have
 \[(q-1) + (2r_S-q+1) + (length \ of \ vertical \ edge) \ge 2r_1+2r_2.\]
 The vertical edge of $P_{q, r_1, q, r_2}$ must have length at least $2r_B$, so the filtration level $(2q, 2q-1)$ supports nontrivial homology.  We claim that this implies $\widehat{HFL}(P_{q_1, r_1, q_2, r_2}, (q_1+q_2, q_1+q_2-1)) \ne 0$ for all $q_i, r_i$.   

      The chain complex $\widehat{CFL}(P_{q_1, r_1, q_2, r_2}, (q_1+q_2, q_1+q_2-1)) = C_1$ has four generators in two gradings.  Each of the top-graded generators can be connected to each of the lower-graded generators by a non-negative domain disjoint from all the basepoints.  It is difficult, however, to determine whether these correspond to boundary discs; although the homology itself is independent of the analytical input, in some cases the choice of complex structure determines whether or not a given domain has a holomorphic representative.  Rather than trying to analyze this directly, we pick the complex structure so that we know the homology is nontrivial.  Specifically, denote by $C_2$ the rank-four chain complex $\widehat{CFL}(P_{q, r_1, q, r_2}, (2q, 2q-1))$ when $r_S > q-1$.  We used a topological argument to show that $C_2$ has nontrivial homology, but the domains connecting generators in $C_1$ are isotopic to domains connecting generators in $C_2$.  Given a generic complex structure $\mathcal{J}_1$ on $\Sigma$, we pick a second complex structure $\mathcal{J}_2$ so that each domain in $(\Sigma, \mathcal{J}_1)$ has the same analytic properties as the corresponding domain in $(\Sigma, \mathcal{J}_2)$.  Although this fails to clarify whether any particular domain has a holomorphic representative,  the answer is the same for corresponding domains in $C_1$ and $C_2$.  Thus, the chain complexes $C_1$ and $C_2$ have the same homology and $\widehat{HFL}(P_{q_1,r_1, q_2, q_2}, (q_1+q_2, q_1+q_2-1) \ne 0$.
      

      $\widehat{HFL}((P_{q_1, r_1, q_2, r_2}, (q_1+q_2, q_1+q_2-1)) \ne 0$ fixes the length of the right vertical edge of $\mathfrak{H}(\widehat{HFL}(P_{q_1, r_1, q_2, r_2}))$ for all $P_{q_1,r_1,q_2, r_2}$ with $r_1 \ne r_2$, and symmetry of the polytope determines the rest of the sides.  The coordinates stated in the theorem then follow from first replacing the temporary filtration levels adopted in Lemma~\ref{filtrpr} with those implied by Equation~\ref{eq:sym} and then applying Theorem~\ref{thm:polytope}.

It remains to show that the same coordinates hold for the case $r_1 =r_2=r$.  The proof that $\widehat{HFL}(P_{q_1,r_1, q_2, q_2}, (q_1+q_2, q_1+q_2-1)) \ne 0$ shows only that the polytope always includes the point $(q_1+q_2, q_1+q_2-1)$; in particular, it does not give the length of the vertical edges for $P_{q_1,r, q_2, r}$.  The filtration level $(q_1+q_2, 2r+q_1+q_2-1)$ supports two generators, but the pair can be connected via a domain satisfying both the positivity and the null intersection conditions.  Explicitly counting points in the moduli space is prohibitively difficult in this case, so we appeal to a new Heegaard diagram.

Beginning with a Heegaard diagram of the type shown in Figure~\ref{fig:P2213} for $P_{q_1, r, q_2, r}$, handleslide $\beta_2$ across $\beta_3$, and then handleslide $\beta_3$ across $\beta_1$.  (See Figure~\ref{fig:HDnew} for an example.)  This creates intersections between previously disjoint curves; in addition to the intersection points of types $A, B, C, a, b, $ and $c$, we have $X_i \ (\alpha_1 \cap \beta_3) $ and $Y_i \ (\alpha_2 \cap \beta_2)$.  Altogether there are four types of generators: $\{ abc, ABC, aAY, cBX \}$.  As before, we assign temporary coordinates consistent with the relative filtration levels of the various intersection points.  (In fact, these are chosen to agree with those in the first diagram, but this is not necessary.  See Figure~\ref{fig:newfiltr}.)  The filtration support is contained in a strip of width $q_1+q_2$, so the right vertical edge of the Heegaard Floer polytope must correspond to the support of the homology on the line $y=q_1+q_2$ in the new coordinates.  
\begin{figure}[h]
Ê \includegraphics[width=\textwidth]{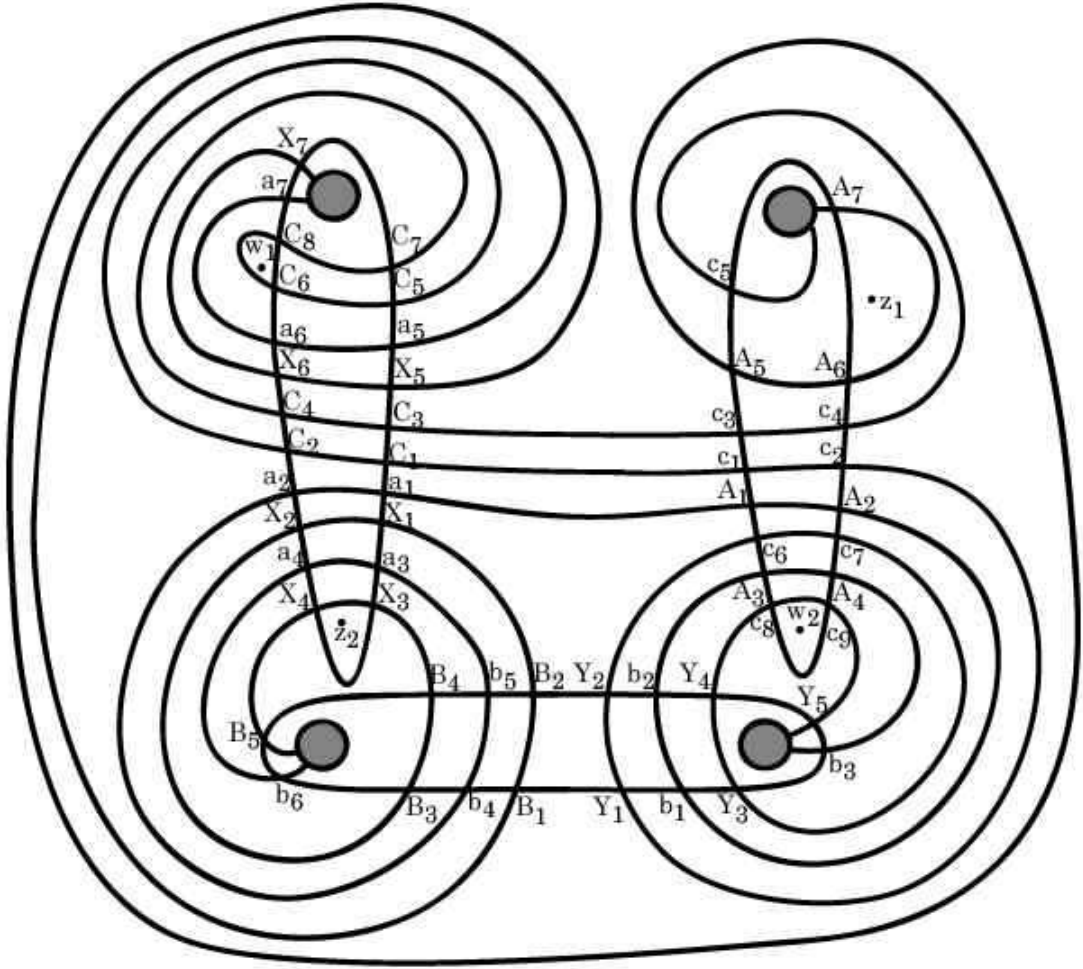}
\caption{A Heegaard diagram for $P_{2,2}$. This diagram is related to the one shown in Figure~\ref{fig:P2213} by two handleslides.}\label{fig:HDnew}
\end{figure} 
As in the first diagram, the filtration level $(q_1 + q_2, q_1+q_2-2)$ has two generators which are connected by a domain supporting a holomorphic disc.  The filtration level $(q_1 + q_2, q_1+q_2-1)$ has four generators, but in contrast to the this filtration level in the first case, we can directly analyze the homology of this complex to show $\widehat{HFL}(P_{q_1, r, q_2, r}, (q_1 + q_2, q_1+q_2-1)) \cong {\Z}_2^2$.  Specifically, there is an obvious simply-connected domain connecting two of the generators, and the techniques of the proof of Lemma~\ref{lem:nonvan} preclude boundary discs between sufficiently many of the others.  Similarly,  $\widehat{HFL}(P_{q_1, r, q_2, r}, (q_1 + q_2, q_1+q_2 +2r-1)) \cong {\Z}_2^2$; there is a unique pair of generators in this filtration level, and no domain connecting them satisfies both the positivity and null intersection conditions.  This shows that the vertical edge has length $2r_B$ for all values of $r_i$ and completes the proof.
\end{proof}
 \begin{figure}
  \includegraphics[width=\textwidth]{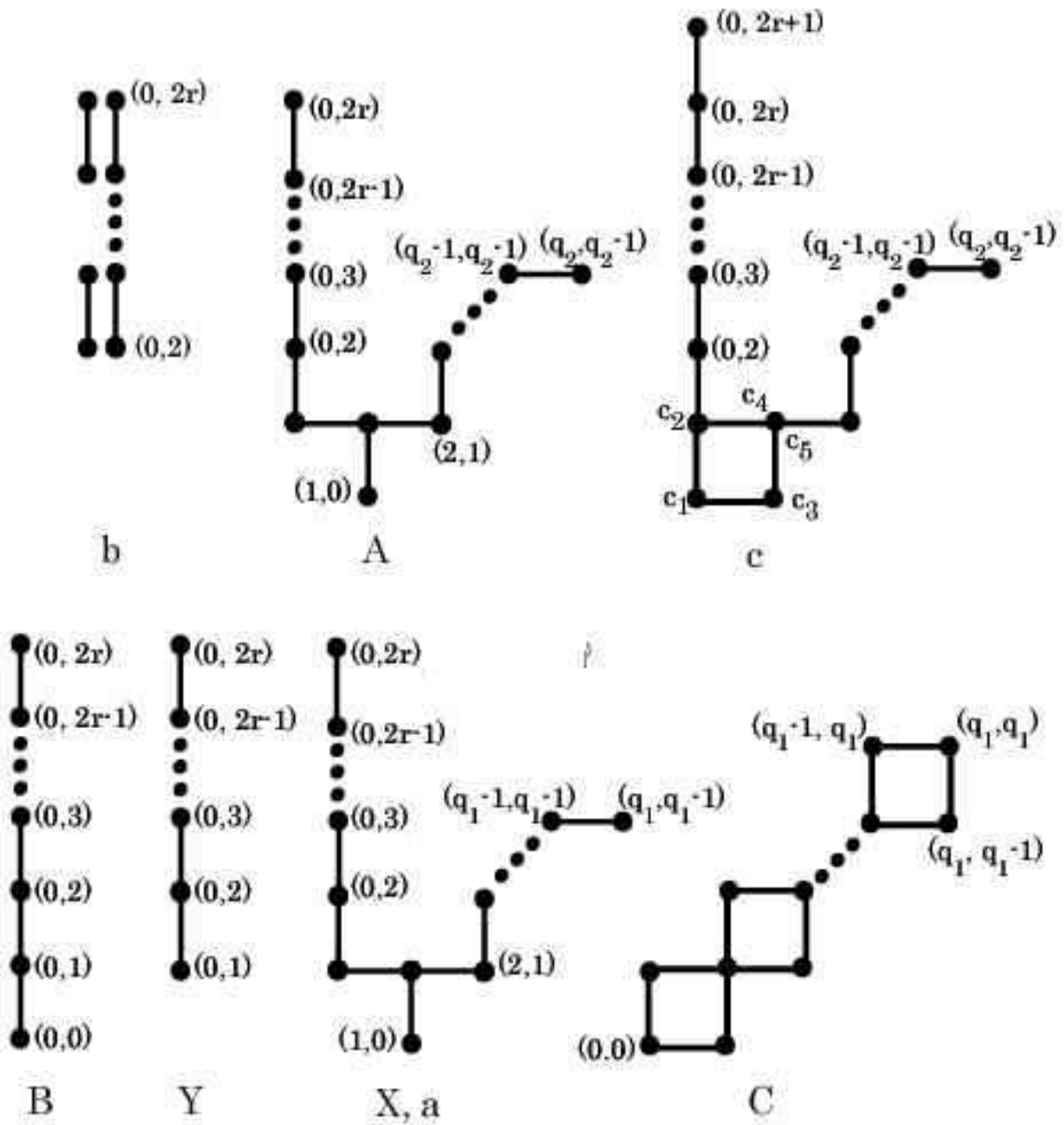}
  \caption{Filtration data for the eight types of intersection points in the type of Heegaard diagram shown in Figure~\ref{fig:HDnew}.}\label{fig:newfiltr}
\end{figure}
The Alexander polynomial of $P_{q_1, r_1, q_2, r_2}$ provides substantial, although not complete, information about the  dual Thurston polytope when $r_1 \ne r_2$.  Because the points $(0, 2r_1+2r_2+1)$ and $(q_1+q_2, q_1+q_2-1)$ have pairs of generators with grading difference one, they do not appear in the Newton polytope.  However, the other extreme points of $S$, where the rank of $\widehat{HFL}(P_{q_1,r_1, q_2,r_2}, \textbf{h})$ is one, correspond to terms in the Alexander polynomial with nonzero coefficients.  In contrast, the Alexander polynomial reveals nothing about  $B_T^*(P_{q,r})$.  Theorem 1 shows that despite this striking difference, minimal spanning surfaces in $P_{q,r}$ behave similarly to those in the general four-parameter family when $2r_S \le q_S-1$: changing a single $r$-crossing changes the complexity of $F_K$ by two regardless of the twist numbers of the other columns.

\end{document}